\documentstyle{amsart}
\input{bull-art}
\newtheorem{theorem}{Theorem}

\newtheorem{proposition}{Proposition}
\newtheorem{corollary}{Corollary}

\theoremstyle{definition}
  \newtheorem{dfn}{\unskip}

\newtheorem{ex}{Example}

\newcommand{\CL}{\operatorname{cl}}
\newcommand{\STAR}{\operatorname{star}}
\newcommand{\ID}{\operatorname{id}}
\newcommand{\Max}{\operatorname{Max}}
\newcommand{\LA}{\operatorname{LA}}
\newcommand{\LS}{\operatorname{LS}}
\newcommand{\Int}{\operatorname{int}}

\begin{document}
\def\currentvolume{29}
\def\currentissue{2}
\def\currentyear{1993}
\def\currentmonth{October}
\def\copyrightyear{1993}
\def\currentpages{189-207}

\title[Intersecting Families of Sets]{Intersecting Families of Sets\\
and the Topology of Cones in Economics}
\author{G. Chichilnisky}
\address{1032 IAB, Columbia University, New York, New York
10027}
\email {gc9@@cunixf.cc.columbia.edu}
\thanks{Research support was provided by NSF Grant No. 
92-16028}
\thanks{The results of this paper were presented at a 
conference on Recent 
Developments in Mathematical Economics sponsored by the 
London Mathematical
Society at the University of Essex, spring 1981, and at 
the European Congress
of Mathematics, Paris, July 1992}
\date{April 20, 1993 and, in revised form, May 21, 1993}

\maketitle
\begin{abstract}

Two classical problems in economics, the existence of a 
market equilibrium
and the existence of social choice functions, are 
formalized here by the
properties of a family of cones associated with the 
economy. 
It was recently established that a necessary and
sufficient condition for solving the former is the 
nonempty intersection of
the family of cones, and one such condition for solving 
the latter is the
acyclicity of the unions of its subfamilies. We show an 
unexpected but clear
connection between the two problems by establishing a 
duality property of
the homology groups of the nerve defined by the family of 
cones. In
particular, we prove that the intersection of the family 
of cones is
nonempty if and only if every subfamily has acyclic 
unions, thus
identifying the two conditions that solve the two economic 
problems. In
addition to their applications to economics, the results 
are shown to extend
significantly
several classical theorems, providing unified and simple 
proofs: Helly's
theorem, Caratheodory's representation theorem, the
Knaster-Kuratowski-Marzukiewicz theorem, Brouwer's fixed 
point theorem, and
Leray's theorem on acyclic covers.
\end{abstract}

%%\tableofcontents
\noindent{1.} Introduction

\noindent{2.} Definitions

\noindent{3.} Market equilibrium

\noindent{4.} Social choice functions

\noindent{5.} Duality and intersecting families

\noindent{6.} Extensions of theorems of Helly, 
Caratheodory, KKM, Brouwer, and Leray

\noindent{7.} Market equilibrium and social choice

\noindent{Acknowledgment}

\noindent{References}
\section{Introduction}
A classic problem in economics is the existence{\em \ }of 
a {\em market
equilibrium} (Von Neumann \cite{vonneumann}, Nash 
\cite{nash}). This can be
viewed as a zero of a nonlinear map $\Psi :$ 
$R^N\rightarrow R^N$
representing market excess demand and embodying optimal 
behavior of the
traders (Arrow and Debreu \cite{arrowdebreu}). The zero 
can be located by
homotopy methods (Eaves \cite{eaves}, Hirsch and\ Smale 
\cite{hirsch2}).
Smale \cite{smale, smale2} has reexamined 
an intuitively appealing
dynamical system which is compatible with a field of cones 
of directions of
improvement for the economy. Along its solution paths all 
traders gain and
proceed until no more gains can be attained and an 
equilibrium is reached.
However, unless the economy satisfies strong boundary 
conditions, this
process may not converge and the market equilibrium may 
fail to exist.

Another classic problem in economics is the existence of 
{\em social choice
functions,} (Arrow \cite{arrow2}). These can be viewed  
(Chichilnisky \cite{chichilnisky}) as maps
which assign to each vector of individual preferences a 
social preference, $
\Phi :P^k\rightarrow $ $P,$ where $P$ is the space of 
preferences and $k$ is
the number of individuals. $\Phi $ must satisfy certain 
properties which
derive from ethical considerations such as symmetry, an 
equal treatment
condition. The problem has a clear
topological structure. A map $\Phi $ exists for a given 
$k$ only when a
certain topological obstruction disappears. It exists for 
all $k$ if and
only if the space $P$ is topologically trivial 
(Chichilnisky and Heal \cite
{chichilniskyheal}). In general, the space $P$ is infinite 
dimensional and has
nontrivial homology, so a social choice rule may fail to 
exist \cite
{chichilnisky, Chichilnisky5}.

Both problems are fundamental to the organization of 
society. Their
solutions model social agreements about how to allocate 
the resources of the
economy among competing individuals, the market solution 
providing an
allocation which is {\em efficient} (Arrow \cite{arrow}) 
and the social
choice solution one which satisfies certain ethical 
properties. The
solutions represent different types of ``social contracts''.

While these two problems appear to be quite different and 
have been
considered separately until now, we show that, in a 
well-defined sense, they
are the same. We provide here a topological formulation of 
these problems
which allows us to identify each with apparently different 
properties of a
family of cones which is naturally associated with the 
economy. 
It was recently shown that the
existence of a competitive equilibrium requires the family 
of cones to
intersect; the existence of social choice functions 
requires that all
subfamilies have acyclic unions (Chichilnisky [12, 13]). 
Looking at the problem in its simplest and
most general form, we obtain a topological 
characterization of a family of
finitely many sets in a general topological space that 
is necessary and sufficient for the
family to have a nonempty 
intersection\footnote{This result was first established 
in Chichilnisky [9].}. 
One main result is that an acyclic or convex family 
has nonempty intersection
if and only if every subfamily has acyclic union (Theorem 
6 and Corollary
2), but the results extend to nonacyclic, nonconvex 
families as well
(Theorems 9 and 10).
As a by-product, we establish the identity between the two 
classical
problems in economics, namely, the existence of a social 
choice function and of a
competitive equilibrium (Theorem 11).

The topology of our family of cones contains crucial 
information about the
economy. The homology of its nerve defines a topological 
invariant for the
economy which provides answers to global problems such as, 
for example,
whether a market equilibrium exists (Theorems 1, Corollary 
2, and Theorem
11). Furthermore, this invariant\ allows us to decide 
whether every
subeconomy has a competitive equilibrium (Theorem 11(b)). 
The homology of
this nerve also contains information about the global 
convergence of the
classic price adjustment process in Smale \cite{smale, 
smale2} (see
Chichilnisky \cite{Chichilnisky11})---it determines 
whether this process
converges.

The homology of the nerve of a family of sets also 
provides valuable 
information in a number of other applications in fields 
other than economics,
which appear as additional by-products of the results in 
this paper. These
include substantial extensions and unified proofs for 
classical
theorems which have until now been considered disparate: 
Helly's theorem on
$n+k$ convex sets in $R^n,$ $k>1$ (\cite{helly, helly2, 
alexandroff}), which is used extensively in game theory, 
for example, 
Guesnerie and Oudu \cite{guesnerie}; Caratheodory's 
theorem and its relative
the Krein-Milman theorem, both of which are used in 
representation theory to
characterize the extreme elements of the cone of positive 
harmonic functions on
the interior of the disk (Choquet \cite{choquet}); the
Knaster-Kuratowski-Marzukiewicz (KKM) theorem (Berge \cite 
{berge}), which is
frequently used to prove the existence of the core of a 
game (Scarf
\cite{scarf}); the Brouwer fixed point theorem, which is 
the nonretractability
of a cell onto its boundary and is used to prove existence 
of solutions of
simultaneous equations (Hirsch \cite{hirsch}, Arrow and 
Hahn \cite{arrowhahn});
and Leray's theorem on the isomorphism between the 
homology groups of a space
and those of the nerves of an acyclic cover (Leray 
\cite{leray}, Dowker
\cite{dowker}, Cartan \cite{cartan}). These classical 
theorems of Helly,
Caratheodory, Leray, and KKM are extended here to simple 
and regular families
of arbitrary finite cardinality, consisting of sets which 
need not be open nor
acyclic or even connected and which are contained in 
general topological
spaces, including  infinite-dimensional spaces; our 
results generalize also the
Brouwer's fixed point theorem which appears as an 
immediate corollary. In
addition, our topological approach allows us to obtain 
conditions which are
simultaneously necessary and sufficient for nonempty 
intersection 
of a general family of sets (Chichilnisky
\cite {chichilnisky10}), a result which we find here very 
useful and which was
not available before.

Here is a summary of the paper. In \S\S 2--4 we set out the
context and describe the problems of existence of a market 
equilibrium and
of a social choice function. A necessary and sufficient 
condition for the
existence of a market equilibrium---called {\em limited 
arbitrage\/}---is
defined as the nonempty intersection of a family of cones. 
A necessary and
sufficient condition for the existence of social choice 
functions---called 
{\em limited social diversity\/}---is defined as the 
acyclicity of the unions
of subfamilies of the same family of cones. Our task is to 
prove that the
two conditions are in fact identical. This identity 
(Theorem 11) 
is a corollary of the
results in \S 5.

Section 5 studies the problem in a general form. First we 
prove a duality
result which relates the reduced singular homology groups 
of the union and
the intersection of a subfamily in dimensions which are 
complementary with
respect to its cardinality \cite{chichilnisky10}. This 
analysis is used to
prove that all subfamilies up to a certain cardinality 
have acyclic unions
if and only if they have acyclic intersections. Then we 
establish that the
whole family has a nonempty acyclic intersection if and 
only if all the
reduced homology groups of the union of its subfamilies up 
to a certain
cardinality vanish.

We further extend the results to families of sets which 
need not be open,
acyclic, or even connected in order to obtain a condition 
for the nonempty
intersection of the family, whether or not this 
intersection is acyclic. The
results thus provide a topological characterization of 
families of sets
which have a nonempty intersection. In particular, this 
characterization
shows that a convex family has a nonempty intersection if 
and only if all
its subfamilies have acyclic unions. Therefore, {\em 
limited arbitrage} is
identical to {\em limited diversity,} and the problems of 
existence of a
competitive equilibrium and of social choice functions are 
the same.

Sections 6 and 7 apply the results in \S 5 to extend a 
number of
classical theorems and to provide simple, unified proofs 
to such disparate
results as Helly's theorem, Caratheodory's representation 
theorem, the
Knaster-Kuratowski-Marzukiewicz theorem, Brouwer's fixed 
point theorem, and
Leray's theorem on acyclic covers. Our extensions of these 
classical results
include families of sets in arbitrary topological spaces 
to which the earlier
results do not apply, sets 
which need not be open, convex, acyclic, or even
connected.  The families may, in addition, be of arbitrary 
finite cardinality.
Section 7 establishes the identity between the problem of
existence of a competitive equilibrium and the problem of 
existence of
social choice functions.

\section{Definitions}

We consider collections of finitely many sets in a 
topological space $X$,
denoted $\{U_\alpha \}_{\alpha \in S},$ with set of 
indices $S$. Such a
collection is called a {\it cover} of $X$ when 
$X=\bigcup_{\alpha \in
S}U_\alpha $; it is an {\em open cover} when each set is 
open in $X$. The
term {\em family} will be used to describe a collection of 
finitely many
sets $\{U_\alpha \}_{\alpha \in S}$ in $X$ whose union 
$\bigcup_{\alpha \in
S}U_\alpha $ may or may not cover X. An {\em open family} 
in $X$ is a family
consisting of sets which are open in $X.$ A subset of 
indices in $S$ will be
indicated by $\theta \subset S$; each subset $\theta 
\subset S$ defines a {\em 
subfamily}{\it \ $\{U_\alpha \}_{\alpha \in \theta }$} of 
the {\em family}
{\it \ $\{U_\alpha \}_{\alpha \in S}$}. We shall use the 
notation $U_\theta $
for the intersection of the subfamily indexed by $\theta 
$, $U_\theta
=\bigcap_{\alpha \in \theta }U_\alpha $, and $U^\theta $ 
for its union $
U^\theta =\bigcup_{\alpha \in \theta }U_\alpha $.

$H_{*}$ will be used to denote {\em reduced} singular 
homology, and $H_q(Y)$
to denote the $q$-singular reduced homology group of the 
space $Y;${\em \
reduced} singular homology is defined by replacing the 
usual chain complex
$$
\cdots C_2\rightarrow C_1\rightarrow C_0\rightarrow 0 
$$
by
$$
\cdots C_2\rightarrow C_1\rightarrow C_0\rightarrow 
Z\rightarrow 0, 
$$
where $Z$ are the integers and $C_0\rightarrow Z$ takes 
each 0-simplex to
1. The corresponding reduced singular homology groups 
denoted $H_q(Y)$ are
defined for all $q\geq -1$. The standard 0-singular 
homology of $Y$ is the
direct sum $H_0(Y)\oplus Z$. Note that with this notation 
if $Y$ is a 
nonempty connected space, then $H_0(Y)=0$ and 
$H_{-1}(Y)=0$; and if $Y$ has two
connected components, then $H_0(Y)=Z$. If $Y$ is empty, 
$H_0(Y)=0$ and $
H_{-1}(Y)=Z$. It is immediate that with this definition 
the Mayer-Vietoris
sequence (Spanier \cite[\S 6, Chapter 4]{spanier}) 
extended to {\em reduced
} singular homology 
$$ 
\cdots H_{q+1}(A\cap B)\rightarrow H_{q+1}(A)\oplus H_{q+
1}(B)\rightarrow
H_{q+1}(A\cup B)\rightarrow H_q(A\cap B)\rightarrow \cdots  
$$
is exact.

We say that a space $Y$ is {\it acyclic} if and only if 
$H_{*}(Y)=0$. Since by
definition the space $Y$ is nonempty 
if and only if $H_{-1}(Y)=0$, in our notation $Y$
is called acyclic when $Y$ is not empty and is acyclic in 
the standard
singular homology. When the space $X$ is contained in a 
linear space, a
family is called {\em convex} if it consists of convex 
sets. A family $
\{U_\alpha \}_{\alpha \in S}$ is called {\it acyclic} if, 
for all $\theta
\subset S$, the set $U_\theta $ is either empty or acyclic.

\begin{dfn} 
For any $k\geq 0$  we say that the family $\{U_\alpha 
\}_{\alpha \in
S}$ satisfies {\bf condition $A_k$} if the intersection 
$U_\theta $ is
acyclic for every $\theta \subset S$ having at most $k+1$ 
elements.
\end{dfn}

\begin{dfn}
For any $k\geq 0$ we say that the family $\{U_\alpha 
\}_{\alpha \in
S}$ satisfies {\bf condition $B_k$} if the union $U^\theta 
$ is acyclic for
every $\theta $ having at most $k+1$ elements.
\end{dfn}

If $X\subset R^n$, then the family $\{U_\alpha \}_{\alpha 
\in S}$ is called
a {\em family in} $R^n$ and is called {\it a }{\em family 
of $k$
sets }if $S$ has cardinality $k$.

\begin{dfn}
If $X$ is a simplicial complex with set of vertices $S$, 
then a {\em 
simple cover} of $X$ is an {\em open cover} $\{U_\alpha 
\}_{\alpha \in S}$
of $X$ satisfying $\CL(U_\alpha )\subset \STAR(\alpha )$ 
for all $\alpha \in S$
, where $\CL(Y)$ is the closure of $Y$ and $\STAR(\alpha 
)$ is the interior
in $X$ of the union of all closed simplices in $X$ having 
$\alpha $ as a
vertex. 
\end{dfn}

The sets in a simple family need not be convex nor acyclic 
 or even
connected. A subcomplex $L$ of a simplicial complex $K$ is 
a subset of $K$
(that is, if $s\in L\Rightarrow s\in K);$ a subcomplex $L$ 
is called {\em 
full} if each simplex of $K$ having all its vertices in 
$L$ belongs to $L$
(Spanier \cite{spanier}). The symbol $[\alpha ]_{\alpha 
\in \theta }$ denotes
the {\em full subcomplex} of $X$ with set of vertices 
$\{\alpha \}_{\alpha
\in \theta }.$

\begin{dfn}
A cover of the simplicial complex $X$ by finitely many 
{\em closed}
sets $\{C_\alpha \}_{\alpha \in S}$ is called{\em \ 
regular} if $\forall
\theta \subset S$, $[\alpha ]_{\alpha \in \theta }\subset 
\bigcup_{\alpha
\in \theta }C_\alpha $.
\end{dfn}

A regular cover $\{C_\alpha \}_{\alpha \in S}$ of a 
simplicial complex $X$
therefore satisfies: for every subset $\theta \subset S$ 
and every simplex $
\Delta $ of $X$ whose vertices lie in $\bigcup_{\alpha \in 
\theta }C_\alpha
, $ we have $\Delta \subset \bigcup_{\alpha \in \theta 
}C_\alpha .$ The sets
in a regular cover need not be convex, acyclic, or even 
connected.

Given a set $X$ and a collection$\{U_\alpha \}_{\alpha \in 
S}$ of subsets of 
$X$, the {\em nerve of }$\{U_\alpha \}_{\alpha \in S}$ is 
the simplicial
complex having as vertices the nonempty elements of 
$\{U_\alpha \}_{\alpha
\in S}$ and whose simplexes are finite nonempty subsets of 
$\{U_\alpha
\}_{\alpha \in S}$ with nonempty intersection (Spanier 
\cite{spanier}).

\section{Market equilibrium}

\subsection[]{A market economy}

A market economy is described by its {\em goods} and its 
{\em traders}.
There are $n$ $>1$ goods and $H>1$ traders. Traders derive 
utility from
vectors (called {\em trades} or{\em \ bundles of goods}) 
in $R^n$, which is
called the {\em consumption }or{\em \ trade space}$.$ Each 
trader is
identified by a vector describing his/her initial {\em 
endowments} of goods $
\Omega _h\in R^n-\{0\}$ and by a real-valued smooth 
$(C^2)$ function $
u_i:R^n\rightarrow R$ which describes the {\em utility} 
derived from the
different consumption vectors. The {\em space of 
allocations} is $R^{nH};$
its elements describe the assignment of one consumption 
vector in $R^n$ for
each trader. The utilities $u_i$ are increasing: $\forall 
x,y\in R^n,$ if $
x\geq y$, then $u_i(x)\geq u_i(y),$ and $\exists \in 
>0:Du_i(x)>\in$, 
where $Du_i(x)$ is the gradient vector of $u_i$ at $x.$ If 
for some $r\in R$
the set $u_i^{-1}(r,\infty )$ is not bounded below in 
$R^n,$ then we assume
that the set of directions of gradients of the 
corresponding hypersurface, $
\{v=Du_i(x)/\Vert Du_i(x)\Vert :$ $u_i(x)=r\}$, is closed 
in $R^n.$ This
assumption is to control the behavior at infinity of the 
leaves of the
foliation of $R^n$ induced by the hypersurfaces of the 
function $u_i$;
geometrically, one rules out ``asymptotic directions'' for 
the gradients on
those hypersurfaces which are not bounded below. A {\em 
market economy} $E$
is therefore defined by its trade space and its traders: 
$E=\{R^n,\Omega
_i,u_i,i=1,\ldots, H\}.$

\subsection[]{Market equilibrium}

Our next tasks are to motivate and then to define the 
notion of a {\em 
competitive equilibrium for the market }$E${\em .} A 
competitive equilibrium
represents a rest point of the trading activity of the 
economy $E$. Trading
requires {\em prices. }A price is a rule which assigns a 
real number called 
{\em value} to  each bundle of goods in a way that depends 
linearly on the
bundles. Therefore, prices are vectors in the dual space 
of the space of
trades, $R^n.$ Each price $p\in R^n$ determines the {\em 
budget set} of a
trader $B(p,\Omega _i)$ consisting of those trades which 
are affordable at
the traders' initial endowment $\Omega _i.$ Therefore, 
$B(p,\Omega _i)=\{x\in
R^n:\langle p,x\rangle=\langle 
p,\Omega _i\rangle\},$ where $\langle.,.\rangle$ is the 
inner product in $R^n.$
Traders trade within their budgets in order to increase, 
ideally to
optimize, their utility.

Trading comes to a rest when a price $p^{*}\in R^n$ is 
found at which the
corresponding set of all optimal trades 
$\{x_i^{*}\}_{i=1,\ldots,H}$ is
compatible with the resources of the economy, i.e., the 
{\em supply} of each
of the $n$ goods equals the {\em demand}. A {\em 
competitive equilibrium} of
the market economy $E$ is therefore defined as a vector of 
prices and of
trades, $(p^{*},x_1^{*}\cdots x_H^{*})\in R^n\times 
R^{nH},$ satisfying the
following conditions: 
%\begin{equation}
\begin{gather}
\label{2}
u_i(x_i^{*})=\operatornamewithlimits
\Max_{x_i\in B(p^{*},\Omega _i)}u_i(x_i) \\
\hfill \text{for }B(p^{*},\Omega _i)=\{x\in R^n:\langle
p^{*},x\rangle=\langle p^{*},\Omega _i\rangle\} \nonumber
\end{gather}
%\end{equation}
and 
%\begin{equation}
\begin{gather}
\label{1}\sum_{i=1}^H(x_i^{*}-\Omega _i)=0\in R^n. 
\end{gather}
%\end{equation}

The vector $x_i^{*}(p^{*})$ is the {\em demand} of trader 
$i$ at prices $
p^{*};$ a solution $x_i(p)$ to problem (\ref{2}) for all 
$p\in R^n$ is the 
{\em demand function} $x_i(p):R^n\rightarrow R^n$ of 
trader $i.$ $
ED(p)=\sum_{i=1}^H(x_i(p)-\Omega _i)$ is the {\em 
aggregate excess demand
function}$
\footnote{The demand and the aggregate excess demand 
functions may not be well defined for some prices which 
are not equilibrium prices.}
$ of the economy{\em \ }$E.$ Condition (2) means that at the
equilibrium allocation all markets clear, i.e., total 
demand for each good
equals total supply, and therefore $ED(p^{*})=0.$

\subsection[]{Market cones}

Consider a market economy $E=\{R^n,\Omega _i,u_i,i=1,  
\ldots,H\}.$ The {\em 
asymptotic preferred cone} $A_i$ is the cone of all 
directions which
intersect every hypersurface of $u_i$ of values exceeding 
$u_i(\Omega _i):$ 
\begin{equation}
\label{a}A_i=\{v\in R^n:\sup _{\lambda \in (0,\infty 
)}\text{ }u_i(\Omega
_i+\lambda v)=\sup _{x\in R^N}u_i(x)\}. 
\end{equation}

The {\em market cone} $D_i$ is 
\begin{equation}
\label{t}D_i=\{p\in R^n:\forall v\in A_i,\langle 
p,v\rangle >0\}. 
\end{equation}
If the utility $u_i$ is a concave function, then both 
cones $A_i$ and $D_i$
are open convex sets, which we now assume. The condition 
of {\em limited
arbitrage} $(\LA)$ is that all {\em market cones} in 
(\ref{t}) intersect:

\begin{equation}
\label{la} (\LA)\quad 
\bigcap_{i=1}^HD_i\neq \emptyset . 
\end{equation}

This means that there exists a price $p\in R^n$ at which 
only limited
increases in utility can be achieved by all traders from 
trades which are
affordable from their initial endowments.

The following has been established:

\begin{theorem}
Limited arbitrage \RM{(\ref{la})} is necessary and 
sufficient for the existence
of a competitive equilibrium in the market $E.$
\end{theorem}

For a proof see Chichilnisky \cite{Chichilnisky4}.

The condition for existence of a competitive equilibrium 
is therefore the
nonempty intersection (\ref{la}) of a family of cones in 
$R^n$ which are
naturally associated with the economy $E,$ namely, of the 
family of market
cones $\{D_i\}_{i=1,\ldots,H}$ defined in (\ref{t}). The 
market cones $
\{D_i\}_{i=1,\ldots,H}$ contain global information about 
the economy, since they
establish directions of utility increases along which all 
utility levels are
eventually reached. As established in Theorem 1, the 
market cones $
\{D_i\}_{i=1,\ldots,H}$ determine whether or not the 
market has a competitive
equilibrium. They also determine whether or not the 
dynamical process 
revisited in 
[34, 35] converges globally; it converges if and only if 
limited arbitrage
holds, i.e., if and only if the family of cones has 
nonempty intersection
(see Chichilnisky \cite{Chichilnisky11}).

The family of market cones $\{D_i\}_{i=1,\ldots,H}$ also 
contains information
about the existence of social choice functions. In the 
next section we shall
see that a condition for existence of a social choice 
function is that every
subfamily of the family of market cones, 
$\{D_i\}_{i=1,\ldots, H},$ has an 
acyclic union. 

\section{Social choice functions}

\subsection[]{Individual and social preferences}

In this section we consider a connected and simply 
connected $CW$ complex $P$
(Spanier \cite{spanier}) representing a space of 
preferences on $R^n$. The
explicit cell structure on $P$ is not needed, only the 
general topological
properties of $CW$ complexes. For example, $P$ could be a 
polyhedron or a
smooth manifold. $P^k$ denotes the product of $P$ with 
itself $k$ times, $
P^k=\stackrel{k\text{ times}}{\overbrace{P\times 
\cdots \times P}}$, and $\Delta P$ is
the ``diagonal'' of $P^k=\{(p_1\cdots p_k)\in P^k:\forall 
i,j=1,\ldots, k,$ $
p_i=p_j\}. $ Examples of spaces of preferences $P$ are 
provided in \S 7.

\subsection[]{Social choice functions}

A {\em social choice function} for the space of 
preferences $P$ and for $k$
individuals, is a continuous map $\Phi :P^k\rightarrow P$ 
assigning to each
vector of $k$ {\em individual preferences} in $P^k$ a {\em 
social preference}
in $P$ satisfying:

1. $\Phi $ is symmetric; i.e., $\Phi $ is invariant under 
the
action of the group of permutations of $k$ letters acting 
naturally on $P^k.$

This condition means that all $k$ individuals are treated 
equally and is
called {\em anonymity}.

2. The map induced by the restriction of $\Phi $ on $\Delta
(P^k)$ at the homotopy level, 
$
(\Phi \mid \Delta (P^k))_{*}:\pi _j(\Delta 
(P^k))\rightarrow \pi _j(P),\text{
is onto }\forall j$.

This condition arises from several applications 
\cite{Chichilnisky5, 
Chichilnisky6}. For example, it is implied by the {\em 
Pareto condition }
\cite{Chichilnisky2aa},{\em \ }which requires that when 
all individuals
prefer one choice $x$ to another $y$, so does society. It 
is also implied by
the assumption that $\Phi \mid \Delta (P^k)=\ID(\Delta 
(P^k))$; i.e., when $
\Phi $ is restricted to the ``diagonal'' of $P^k,$ $\Delta 
P^k=\{(\rho
_1\cdots \rho _k)\in P^k$ s.t. $\forall i,j,$ $\rho 
_i=\rho _j\},$ it is the
identity map$.$ This latter condition means that when all 
individuals have
the same preference, society adopts that common 
preference, and it is called 
{\em respect of unanimity }\cite{chichilnisky}.

An allocation is an assignment of a bundle of goods in 
$R^n$ to each
trader, and the space of allocations is $R^{nH}$. Each 
trader has a
preference over allocations. A {\em smooth preference} 
over the space of
allocations $R^{nH}$ is a smooth ($C^2)$ unit vector field 
$\rho
:R^{nH}\rightarrow S^{nH-1}$ satisfying: $\exists 
u:R^{nH}\rightarrow R$
with $\forall x\in R^{nH},$ $\rho (x)=\lambda (x)Du(x)$ 
for some $\lambda
(x)>0$ (Debreu \cite{debreu}). The space of all smooth 
preferences on
allocations in $R^{nH}$ is denoted $\Gamma (R^{nH}).$ {\em 
The space of
preferences $P(E_\theta )$ similar to those of a subset 
$\theta \subset
\{1,\ldots, H\}$} {\em of traders in $E$} is
\begin{equation}
\label{u}P(E_\theta )=\left\{\rho \in \Gamma 
(R^{nH}):\text{ }\forall x\in R^{nH} 
\text{ and }\forall j,\text{ }\rho ^j(x)\in \bigcup_{i\in 
\theta }D_i\right\}, 
\end{equation}
where $\rho ^j(x)$ is the projection of $\rho (x)$ on the 
$j$\<th
copy of $R^n$ in the product space $R^{nH}.$ The 
interpretation is that $
P(E_\theta )$ consists of all preferences which are 
similar to those of some
trader $i\in \theta $ in some position $j$ in the sense 
that they
increase in the directions of large utility increases for 
$i$ in position $
j $ and only in those directions. This is discussed 
further in \S 7.
Note that the notion of similarity of preferences depends 
on the same family
of {\em market cones} $\{D_i\}_{i=1,\ldots,H}$ defined in 
equation (4) in \S 3.  

\subsection[]{Social choice and the topology of preferences}

In its most general form the problem of existence of 
social choice
functions has no solution; for the space $\Gamma =\Gamma 
(R^m)$ of all
smooth preferences on $R^m,$ $m>2:$

\begin{theorem}
There exists no map $\Phi :\Gamma ^k\rightarrow \Gamma $ 
satisfying 
\RM{4.2.1} and \RM{4.2.2}
$\forall k\geq 1.$
\end{theorem}

A proof is in Chichilnisky \cite{chichilnisky, 
Chichilnisky2aa}.

A natural question is what spaces of preferences $P$ admit 
a social choice
function. The following is known:

\begin{theorem}
There exists a social choice map $\Phi :P^k\rightarrow P$ 
satisfying \RM{4.2.1}
and \RM{4.2.2} $\forall k\geq 1,$ if and only if $P$ is 
acyclic.
\end{theorem}

This was proved in Chichilnisky \cite{chichilnisky} and 
Chichilnisky and
Heal \cite{chichilniskyheal}.

When a social choice function $\Phi :P^k\rightarrow P$ 
exists, then by
Whitehead's theorem (Spanier [36]) $P$ is contractible, 
since
the space $P$ is acyclic and by assumption $\pi _1(P)=0$. 
%by Whitehead's theorem (Spanier \cite{spanier}) $P$ is 
%contractible. 
Therefore, there exists
a continuous deformation of the space of preferences $P$ 
into one
preference.  For this reason, in this context the 
acyclicity of a space of
preferences establishes a{\em \ limit on social diversity} 
(Heal \cite{heal}).
For any given subset $\theta $ of traders in $E,$ $\theta 
\subset
\{1,\ldots, H\},$ a social choice function $\Phi :\ 
(P(E_\theta ))^k\rightarrow
P(E_\theta )$ exists satisfying the required conditions 
$\forall k>1$ if
and only if $P(E_\theta )$ is acyclic. This in turn means 
that the space of
gradients of the preferences in $P(E_\theta ),$ namely, 
$\bigcup_{i\in \theta
}D_i,$ must be acyclic. We say the market $E$ has {\em 
limited social
diversity} or{\em \ }simply{\em \ limited diversity } 
$(\LS)$, when: 

\begin{equation}
\label{ls}(\LS)\quad 
\forall \text{ }\theta \subset \{1,\ldots, H\},\quad 
\theta \neq
\emptyset \!\Rightarrow \text{ }\bigcup_{i\in \theta 
}D_i\text{ is acyclic.} 
\end{equation}
A consequence of Theorem 3 is:

\begin{theorem}
There exists a social choice function $\Phi :\ P(E_\theta 
)^k\rightarrow
P(E_\theta )$ satisfying \RM{ 4.2.1} and \RM{ 4.2.2}, 
$\forall \theta \subset
\{1,\ldots, H\}$ and $\forall k\geq 1$, if and only if the 
market $E$ has
limited social diversity $(\LS).$
\end{theorem}

This follows from Chichilnisky \cite{chichilnisky} and 
Chichilnisky and
Heal \cite{chichilniskyheal}.

\subsection[]{Social choice and the nerve of market cones}

For a social choice function $\Phi $ to exist, the union 
of every nonempty
subfamily of market cones $\{U_i\}_{i=1,\ldots, H}$ must 
be acyclic$.$ We saw in
\S 3 that the existence of a competitive equilibrium 
requires the 
nonempty intersection of the same family of market cones, $
\bigcap_{i=1}^H\{D_i\}\neq \emptyset .$ To identify the 
two economic
problems, we must exhibit the connection between two 
properties of the
family of cones. One is that the family has nonempty 
intersection---i.e., 
{\em limited arbitrage} (\ref{la}). The second is that the 
union of every
subfamily is acyclic---i.e., {\em limited diversity} 
(\ref{ls}). This is
achieved in Theorem 11 in \S 7 and motivates the results 
in the
following section. 

\section{Duality and intersecting families}

Having established the importance in economics of the 
topology of the nerve
of the {\em market cones} $\{D_i\}_{i=1,\ldots, H},$ we 
turn now to the
mathematical problem. In their simplest and most general 
form the questions
are: 
when does the family of market cones $\{D_i\}_{i=1,\ldots, 
H}$ have a nonempty
intersection, and how does this relate to the acyclicity 
of the unions of its
subfamilies? The nonempty intersection of this family of 
cones is the
condition of {\em limited arbitrage} (\ref{la}), and the 
acyclicity of the
unions of its (nonempty) subfamilies is the condition of 
{\em limited
diversity} (\ref{ls}). We saw in \S 3 that the former 
(\ref{la}) is
necessary and sufficient for the existence of a market 
equilibrium and in
\S 4 that the latter (\ref{ls}) is necessary and 
sufficient for the
existence of social choice functions. This section will 
establish {inter
alia} that the two mathematical conditions (\ref{la}) and 
(\ref{ls}) are
identical.

Here is a summary of the 
section.\footnote{The results in this 
section were first established in Chichilnisky (1981).}
Theorem 5 proves the equivalence between two topological 
conditions of the
nerve of a family of sets of a general topological space 
$X$\<---these are
conditions $A_k$ and $B_k$ defined in \S 2, the former 
requiring that
all subfamilies with at most $k+1$ elements have acyclic 
intersection and
the latter requiring that all such subfamilies have 
acyclic unions. This
identity is simple and geometrically appealing. It has 
many implications, as
we show below. Because it is close to the foundations of 
homology theory,
there is a subtle point in its proof, which ensures an 
excision property
for singular reduced homology (see, e.g., Spanier 
\cite[p.\ 189]{spanier}) so
that the Mayer-Vietoris sequence for reduced singular 
homology---a 
sequence which is rarely used for families where the sets 
may have empty
intersection---is exact.
A discussion of this exactness for reduced homology for
families which includes empty sets is in \S 2, and the 
excision property
is discussed in this section after condition (6).

The exactness of the Mayer-Vietoris sequence is used in 
our proof of a
duality property of the singular reduced homology of a 
family of sets in
Proposition 1. This proposition establishes a simple 
isomorphism between the
singular reduced homology groups of the union and those of 
the intersection of a 
subfamily in dimensions complementary with its 
cardinality. This duality
property allows us to prove 
the following 
somewhat surprising result in Proposition 2: For
families in $R^N$ the conditions $A_k$ and $B_k$ need only 
be required for
subfamilies with at most $N+1$ sets; they are 
automatically satisfied
otherwise$.$ The geometric implications of these results 
are shown in 
Corollary $1,$ which shows that, if the family is acyclic 
and every
subfamily with at most $N+1$ sets has a nonempty 
intersection, then the
whole family has a nonempty intersection.

Building on this, Theorem 6 gives a
necessary and sufficient condition for the acyclic (and 
therefore nonempty)
intersection of 
every subfamily of a family of finitely many sets in a 
general topological
space $X$; the union of every subfamily must be acyclic. 
Furthermore, if the
family of sets is in $R^N,$ the acyclicity is required 
only for subfamilies
with no more than $N+1$ sets. For acyclic families, 
Corollary 2 gives a
simple, necessary, and sufficient condition for the 
nonempty intersection of
the whole family; particularly,
the family has nonempty intersection if and only if every
subfamily has an acyclic union. This result is just what 
is needed for the
economic applications presented in \S\S 3 and 4, as seen 
in Theorem 11
in \S 7.

So far we have considered families which have either empty 
or acyclic
intersection and have excluded those where the sets have 
nonempty intersection,
but this intersection is not acyclic. In several 
applications, for example, 
for non-convex economies it is necessary to consider the 
condition of
limited arbitrage (\ref{la}) which requires nonempty 
intersection, even
when this intersection fails to be acyclic. Here 
Mayer-Vietoris is no longer
useful, and other arguments are needed. 
The rest of this section extends
the results to families which may have nonacyclic as well 
as nonempty
intersection. This is achieved as follows: Theorems 7 and 
8 establish an
isomorphism between the homology of a space $X$ and that 
of the nerve of a 
{\em simple} and of a {\em regular} cover respectively, as 
defined in
\S 2. These include covers by sets which may be neither 
open, convex,
acyclic, or even connected. Using this isomorphism, 
Theorems 9 and 10 prove
necessary and sufficient conditions for nonempty 
intersection; these are
similar to Theorem 6, but they are valid for simple and 
for regular families
respectively.

Unless otherwise stated, the following results apply to a 
general
topological space $X,$ and the family $\{U_\alpha 
\}_{\alpha \in S}$
satisfies 
\begin{equation}
\label{ex}\bigcup_{\alpha \in S}U_\alpha =\bigcup_{\alpha 
\in S}\left(
\Int_{U^S}\text{ }(U_\alpha )\right) , 
\end{equation}
where $\Int_{U^S}(U_\alpha )\text{ denotes the interior of 
the set }
U_\alpha $ relative to the set $U^S=$ 
$\bigcup_{\beta \in S}(U_\beta )$.
A family satisfying this property (\ref{ex}) is called an 
{\em excisive
family}. Since we can take $X=\bigcup_{\alpha \in 
S}U_\alpha ,$ (\ref{ex})
is a rather general specification. For example, (\ref{ex}) 
is satisfied
when the family consists of sets $U_i\subset X$, each of 
which is open in $X$.
Note, however, that condition (\ref{ex}) does not require 
that the sets $U_i$
be open in $\bigcup_{\alpha \in S}(U_\alpha ).$ In fact, 
(\ref{ex}) is
strictly weaker than the requirement that the sets $U_i$ 
be open in $X$; it
includes, for example, families consisting of two closed 
sets $C_1$ and $C_2$
in $R^n$ with $C_1\subset C_2.$ The role of (\ref{ex}) is 
to ensure the
union and the intersection of any subfamily of $\{U_\alpha 
\}_{\alpha \in S}$
define an excisive couple so that the Mayer-Vietoris 
sequence of reduced
singular theory is exact (see Spanier \cite[Theorems 3, 4 
and 
Corollary 5, pp.\ 188--189]{spanier}). 
An example in \cite[p.\ 188]{spanier} 
exhibits two closed path-connected sets $Y_1$ and $Y_2$ in 
$R^2$ such that $%
Y_1\cup Y_2=R^2$ which do not satisfy (\ref{ex}) and for 
which the
corresponding singular Mayer-Vietoris sequence is not 
exact. Condition (\ref
{ex}) prevents such pathologies.

\begin{theorem}
An excisive family $\{U_\alpha \}_{\alpha \in S}$ in $X$ 
satisfies $A_k$ if
and only if it satisfies $B_k.$
\end{theorem}

\begin{pf} The first step in the proof is to establish the
following duality result:
\renewcommand{\qed}{}\end{pf}

\begin{proposition}
Consider an excisive family of sets in $X,$ $\{U_\alpha 
\}_{\alpha \in S},$
satisfying $A_{k-1}$, for $k\geq 1$. Then if $\theta 
\subset S$ has $k+1$
elements, for all $q\ $ 
\begin{equation}
\label{d}H_q(U^\theta )\simeq H_{q-k}(U_\theta ).\text{ } 
\end{equation}
\end{proposition}

\begin{pf} We proceed by induction. When $k=1$, the family 
has
two sets, and this is the Mayer-Vietoris sequence for 
reduced singular
homology as defined in \S 2. Assume the result is true for 
every family 
$\{U_\alpha \}_{\alpha \in \theta }$ where $\theta $ has 
$k$ elements.
Consider now a family $\{U_\alpha \}_{\alpha \in \tau }$ 
of $k+1$ elements 
satisfying $A_{k-1}$. Define $\theta $ so that $\tau 
=\{0\}\cup \theta $,
and $V_\alpha =U_0\cup U_\alpha $ , $\alpha \in \theta $ . 
The new family $%
\{V_\alpha \}_{\alpha \in \theta }$ has $k$ elements, and 
it satisfies $%
A_{k-2}$ because the family 
$\{U_\alpha \}_{\alpha \in \tau }$ satisfies $A_{k-1}$ and
by Mayer-Vietoris. Then 
\begin{align*}
H_q(U^\tau )&=H_q(V^\theta )=H_{q-(k-1)}(V_\theta )
\quad \text{by the induction hypothesis} \\
&=H_{q-k+1}(U_0\cup [U_1\cap \cdots \cap U_k])\\
&=H_{q-k}(U_0\cap [U_1\cap \cdots \cap
U_k])\quad \text{by Mayer-Vietoris}\\ 
&=H_{q-k}(U_\tau ),
\end{align*}
completing the proof of the proposition. The rest of the 
proof of Theorem 5
follows from Proposition 1 by induction on $k$. 
%%$\diamondsuit $
\end{pf}

\begin{proposition}
Let $\{U_\alpha \}_{\alpha \in S}$ be an excisive family 
in $R^n$ satisfying 
$A_n$ . Then $\{U_\alpha \}$ also satisfies $A_k$ and 
$B_k$ for all $k\geq n$%
. In particular, the intersection of this family is always 
nonempty.
\end{proposition}

\begin{pf} This follows from Theorem 5 and Mayer-Vietoris,
because $H_i(U)=0$ for $i\geq n$ for an open set $U\subset 
R^n$.%
%$\diamondsuit $
\end{pf}

\begin{corollary}
Let $\{U_\alpha \}_{\alpha \in S}$ be an acyclic excisive 
family in $R^n$
with at least $n+1$ elements. If every subfamily with $n+
1$ elements has
nonempty intersection, then the whole family has a 
nonempty intersection.
\end{corollary}

\begin{pf} This follows from Proposition 2 because $A_n$ is
satisfied by acyclicity. %$\diamondsuit $
\end{pf}

\begin{ex} 
The conditions of Proposition 2 and Corollary 1 cannot be
relaxed. In general, 
the family must have finite cardinality. Consider, for 
example, the infinite
family in $R^1\{U_i\}_{i=1,2,\ldots}, 
U_i=(i,\infty ).$ Every subfamily of $%
\{U_i\}_{i=1,2,\ldots}$ has acyclic union, but the whole 
family has empty
intersection. Figure 1 shows that Corollary 1 does not 
hold for nonacyclic
families; each three of these four sets in Figure 1 
intersect, but the whole
family has an empty intersection. Figure 2 also shows that 
Proposition 2 is
not true when $A_n$ is not satisfied. Here $n=2$, and 
$A_2$ is not satisfied
because the union of two of the sets is not acyclic. 
%
%\FRAME{dtbpF}{3.49984in 
%}{3.49984in}{0in}{}{}{blank.eps}{\GRAPHICSPS{%
%\FILENAME{C:/BLANK.EPS}}}
%
\end{ex}

  \begin{figure}
  \noindent % otherwise figures would be vertically stacked
  \begin{minipage}{.4\textwidth}
  \vspace{11.5pc}
  \caption{}
  \end{minipage}
  \hfill
  \begin{minipage}{.4\textwidth}
  \vspace{11.5pc}
  \caption{}
  \end{minipage}
  \end{figure}

\begin{theorem}
Let $\{U_\alpha \}_{\alpha \in S}$ be an excisive family 
of $k\geq 2$ sets.
Then the intersection of every subfamily 
$\bigcap_{\alpha \in \theta }U_\theta, \forall 
\theta\subset S $, is acyclic \RM(and
hence nonempty\/\RM) if and only if the union of every 
subfamily $\bigcup_\theta
U_\alpha $ , $\forall \theta \subset S$, is acyclic\/\RM; 
i.e., the family
satisfies{\bf \ $B_{k-1}$}\nolinebreak. 
If the family $\{U_\alpha \}_{\alpha \in S}$ is
in $R^n$, then its intersection is acyclic if and only if 
its union $%
\bigcup_{\alpha \in S}U_\alpha $ is acyclic and it 
satisfies $B_j$ for $%
j=\min (n,k-2)$.
\end{theorem}

\begin{pf} The first statement follows from Theorem 5. For 
the
second statement, first let $j=k-2$. Assume that 
$\bigcup_{\alpha \in
S}U_\alpha $ is acyclic and $\{U_\alpha \}_{\alpha \in S}$ 
satisfies $%
B_{k-2}$. Then $B_{k-1}$ is satisfied. By Theorem 5 so is 
$A_{k-1}$ so
that the intersection of the family is acyclic and thus 
nonempty.
Reciprocally, if the intersection of the whole family is 
not empty, then $%
A_{k-1}$ is satisfied and by Theorem 5 so is $B_{k-1}$ so 
that the union of
the family is acyclic. Now let $j=n$. By assumption and 
Theorem 5, $A_n$ is
satisfied. By Proposition 2 this implies that the whole 
family has nonempty
intersection and that $A_m$ is satisfied for all $m\geq 
0$. Therefore by
Theorem 5, $B_m$ is satisfied for all $m$, and the 
family's union is
acyclic. %$\diamondsuit $
\end{pf}

\begin{ex}
Figures 2 and 3 show that the conditions of Theorem 6
cannot be relaxed. Figure 2 shows that ``acyclic 
intersection'' cannot be
replaced by ``nonempty intersection''; it depicts two sets 
which do
intersect but have a nonacyclic union. Figure 3 shows that 
Theorem 6 is
not true if we replace ``acyclic union'' by ``contractible 
union'' in its
statement; it depicts two ``comb'' spaces having an 
acyclic (and hence
nonempty) intersection, the point $\{x\}$. The union of 
the two comb spaces
is acyclic, confirming Theorem 6, but it is not 
contractible.
%Example 2: %
%\FRAME{dtbpF}{%
%3.49984in}{3.49984in}{0in}{}{}{blank.eps}{\GRAPHICSPS{%
%\FILENAME{C:/BLANK.EPS}%
%}}
\end{ex}
\begin{corollary}
An acyclic excisive family $\{U_\alpha \}_{\alpha \in S}$ 
has
nonempty intersection if and only if $\forall \theta 
\subset S$, the union
of the subfamily $\{U_i\}_{\alpha \in \theta 
},\bigcup_{\alpha \in
S}U_\alpha ,$ is acyclic.
\end{corollary}

\begin{pf} This follows from Theorem 6 and the definition of
acyclic families. %$.\diamondsuit $
\end{pf} 

  \begin{figure}
  \noindent % otherwise figures would be vertically stacked
  \begin{minipage}{.45\textwidth}
  \vspace{9pc}
  \caption{}
  \end{minipage}
  \hfill
  \begin{minipage}{.45\textwidth}
%  \vspace{9pc}
\vspace{9pc}
  \caption{}
  \end{minipage}
  \end{figure}

\begin{ex}
The conditions of Corollary 2 cannot be relaxed. Figure 4 
depicts
a family of $k=4$ sets in $R^2$ which does not satisfy 
$B_2$ (or $A_2)$
because three of them do not intersect. The union of the 
family is acyclic,
but the intersection is empty.
%
%\FRAME{dtbpF}{3.49993in}{3.49993in}{0in}{}{}{%
%blank.eps}{\GRAPHICSPS{\FILENAME{C:/BLANK.EPS}}}
%

Until now we considered
families which had either empty or acyclic intersection. 
The following
results apply to {\em simple} and {\em regular} families, 
as defined in
\S 2. These may consist, for example, of sets in $R^n$ 
which are
neither open nor acyclic or even connected. The families 
may have
{\em nonacyclic, nonempty intersection}. Mayer-Vietoris is 
not useful in this
context, and we must adopt a different approach.
\end{ex}

If $X$ is a simplicial complex, the expression $X$ $=$ 
nerve $\{U_\alpha
\}_{\alpha \in S}$ is used to indicate that $X$ and nerve 
$\{U_\alpha
\}_{\alpha \in S}$ have the same combinatorial structure.

\begin{theorem}
Let $\{U_\alpha \}_{\alpha \in S}$ be a simple cover of a 
simplicial complex 
$X$ with set of vertices equal to $S$. Then $X$ $=$ 
$\{U_\alpha \}_{\alpha
\in S}$.
\end{theorem}

\begin{pf} The proof follows by induction on the number of 
sets $k$. Let the
set of vertices $S$ consist of $k=2$ elements. Then $X$ is 
either a segment
or a set of two points; assume $X$ is a segment. Consider 
$x\ \in \ \partial
U_1$. Since $x\notin U_1$, $x\ \in \ U_2$. Therefore, 
$\exists y\ \in \
U_1\cap U_2$. Now let $X=\{x_1\}\cup \{x_2\}$. Since 
$U_\alpha \subset
\STAR(\alpha )$, $U_1\cap U_2$ is empty. Consider now the 
following inductive
assumption for a set of vertices $S$ of $k+1$ elements: 
the nerve $%
\{U_\alpha \}=X$, and if the $k$ sets $\{U_\alpha 
\}_{1\leq \alpha \leq k}$
intersect, then $\exists $ a simple family $\{W_\alpha 
\}_{1\leq \alpha \leq
k+1}$ covering $X$ with $W_\alpha \subset U_\alpha \forall 
\alpha $ and an $%
x\ \in \ \partial W_1\cap \cdots \cap \partial W_k$. Now 
let $S$ have $k+2$
sets. Assume $X$ is a $k+1$ simplex. By the inductive 
hypothesis every
subfamily of $k+1$ sets in $\{U_\alpha \}$ intersects, and 
in particular, $%
\exists $ a simple family $\{W_\alpha \}_{1\leq \alpha 
\leq k+1}$ with $x\in
\bigcap_{1\leq \alpha \leq k}\partial W_\alpha $ . Let 
$Z_{k+1}=W_{k+1}-I_x$%
, where $I_x$ is a closed segment in $W_1\cap \cdots \cap 
W_k$, and $x\ \in \
\partial I_x$. Take $Z_{k+1}$ to be an element of the 
simple family $%
\{Z_\alpha \}_{1\leq \alpha \leq k+2}$ defined otherwise 
by $Z_\alpha
=W_\alpha $ for $\alpha \leq k$ and $Z_{k+2}=U_{k+2}$. 
Then $\forall \alpha 
$, $Z_\alpha \subset U_\alpha $, $\{Z_\alpha \}_{1\leq 
\alpha \leq k+2}$
covers $X$, and $x\ \in \ \partial Z_1\cap \cdots \cap 
\partial Z_{k+1}$, so
$x\ \in \ Z_{k+2}$. Therefore, $x\ \in \bigcap_{1\leq 
\alpha \leq 
k+2}Z_\alpha \subset \bigcap_{1\leq \alpha \leq k+
2}U_\alpha \neq \emptyset $.
Finally, if $X$ is not a simplex, $\bigcap_{1\leq \alpha 
\leq k+2}U_\alpha
=\emptyset $, since $U_\alpha \subset \STAR(\alpha )$ for 
all $\alpha \in S$. 
%$\diamondsuit $
\end{pf}

The following result uses the definition of regular covers 
given in \S
2.

\begin{theorem}
Let $\{C_\alpha \}_{\alpha \epsilon S}$ be a regular cover 
of a simplicial
complex $X.$ Then nerve $\{C_\alpha \}_{\alpha \in S}=X$ .
\end{theorem}

\begin{pf} First we prove that Theorem 7 implies that if $%
\{C_\alpha \}_{\alpha \in S}$ is a regular cover of $X$, 
then $%
\bigcap_{\alpha \in S}C_\alpha \neq \emptyset .$ Let 
$D_\alpha =C_\alpha \cap
\STAR(\alpha )$; then $\bigcup_{\alpha \in S}D_\alpha =X$. 
Now by Theorem 7
\begin{equation}
\text{if $\{U_\alpha \}_{\alpha \in S}$ is a simple family 
covering $X$ with $
U_\alpha \supset D_\alpha $ for all $\alpha $, 
$\bigcap_{\alpha \in
S}U_\alpha \neq \emptyset $}. 
\end{equation}
We now use (10) to prove $\bigcap_{\alpha \in
S}D_\alpha \neq \emptyset $ , by induction on $k$. 

Case $k=1$. If $%
\bigcap_{\alpha =1,2}D_\alpha =\emptyset $, then $\exists 
U_1, U_2$ defining a
simple family with $\bigcap_{\alpha =1,2}U_\alpha 
=\emptyset $, contradicting
(10). Now let $S$ have $k+1$ elements: by the inductive 
assumption, $%
\bigcap_{1\leq \alpha \leq k}C_k\neq \emptyset $. If 
$[\bigcap_{1\leq \alpha \leq
k}C_\alpha ]\cap C_{k+1}=\emptyset $, then $\exists $ a 
simple family $\{U_\alpha
\}$ s.t. $[\bigcap_{1\leq \alpha \leq k}U_\alpha ]\cap 
U_{\alpha +1}=\emptyset $,
contradicting (10). 
Thus $\bigcap_{\alpha \in S}D_\alpha \neq \emptyset $ so 
that $%
\bigcap_{\alpha \in S}C_\alpha \neq \emptyset .$ 

Having established the result
for the case where $X$ is a simplex, the rest of the proof 
follows the proof
of Theorem 7 by considering the family defined by the 
complements of the
sets $\{C_\alpha \}_{\alpha \in S}$ in $X$.%$\diamondsuit $
\end{pf}

The two following theorems extend the results of Theorem 6 
to the cases of 
{\em simple} and {\em regular} families as defined in \S 
2; here we are
concerned with the nonempty intersection of the family, 
whether or not this
intersection is acyclic.

\begin{theorem}
Let $\{U_\alpha \}_{\alpha \in S}$ be a simple family of 
$k$ sets, such that
every subfamily with $k-1$ elements has a nonempty 
intersection. Then the
whole family has a nonempty intersection if and only if 
its union $%
\bigcup_{\alpha \in S}U_\alpha $ is acyclic. If $k>n+1$, 
we need to
require only that every family of $n+1$ sets has a 
nonempty intersection.
\end{theorem}

\begin{pf} By assumption the $(k-2)$-skeleton of nerve $
\{U_\alpha \}_{\alpha \in S}$ is the boundary of a $k-1$ 
simplex. Let $
X=\{U_\alpha \}_{\alpha \in S}$. By Theorem 7 nerve 
$\{U_\alpha \}_{\alpha
\in S}=X$. Therefore, all sets in the family $\{U_\alpha 
\}$ intersect if and
only if its union $X=\bigcup_{\alpha \in S}U_\alpha $ is 
acyclic. 
%$\diamondsuit $
\end{pf}

\begin{theorem}
Let $\{C_\alpha \}_{\alpha \in S}$ be a family of $k$ 
closed sets with $%
[\alpha ]_{\alpha \in \sigma }\subset \bigcup_{\alpha \in 
\sigma }C_\alpha $%
, and $\bigcap_{\alpha \in \sigma }C_\alpha \neq \emptyset 
$ for every subset $%
\sigma $ of $S$ with $k-1$ elements. Then  
$\bigcap_{\alpha \in \sigma
}C_\alpha \neq \emptyset $ if and only if $\bigcup_{\alpha 
\in \sigma }C_\alpha $
is acyclic. If $k>n+1$, we need to require only that every 
family of $n+1$
sets has a nonempty intersection. 
\end{theorem}

\begin{pf} This follows from the proof of Theorem 9,
%substituting in its proof 
replacing Theorem 7 in the proof by Theorem 8. 
%%$\diamondsuit $
\end{pf}

\section[Extensions of theorems of Helly, Caratheodory, 
KKM, Brouwer, and
Leray]{Extensions of theorems of Helly, Caratheodory, 
KKM,\protect\\ Brouwer,
and Leray}

The question of when sets intersect was studied in the 
classic theorems of
Helly \cite{helly, helly2} and of 
Knaster-Kuratowski-Marzukiewicz in 
\cite{berge}. They provided conditions which are 
sufficient for a
family of sets in $R^n$ to have a nonempty intersection, 
but their results
are restricted to families with $n+1$ or more sets in the 
case of Helly's
theorem and to families with exactly $n+1$ sets in the 
case of KKM's
theorem, in both cases having either a convex structure or 
other particular
characteristics. These two results are quite specific to 
the problems they
study and appear to be different from each other. However, 
the problem of
nonempty intersection in its general form has a clear 
geometrical
structure and can be dealt with by using topological 
tools. We showed in
\S 5 that no restrictions on the number of sets is 
required, nor is
convexity, acyclicity, or even connectedness of the sets. 
Furthermore, the
families need not be in $R^n$ or in any linear space. Once 
this is
understood, the two classic results of Helly and KKM 
appear as special cases
of our results. Brouwer's theorem is also a special case 
of our results,
since it is known to be implied by the KKM theorem, as is 
Caratheodoty's
theorem, which follows from the Helly's theorem.

Helly's theorem is connected here to the Brouwer's fixed 
point theorem and
to an extension provided here of Leray's theorem on 
acyclic covers. Our
extension of Leray's theorem (Leray \cite{leray}, Dowker 
\cite{dowker},
Cartan \cite{cartan}) is in Theorems 7 and 8 of \S 5; 
while Leray's
theorem applies to acyclic covers and proves the 
isomorphism of the
homology of the nerve of the cover and that of the union 
of the family, our
Theorems 7 and 8 significantly extend
this result for covers consisting of sets which may
not be acyclic nor open or even connected.

This section therefore exhibits how the results in \S 5 
extend and
unify several classical theorems. Proposition 2 in \S 5 
extends
Helly's striking theorem on the nonempty intersection of 
families in $R^n$
having more than $n+1$ sets (Helly \cite{helly}, 
Alexandroff and Hopf \cite
{alexandroff}) to possibly nonconvex and nonacyclic 
families with any
number of sets in a general topological space $X$. 
Corollary 3 below is
Helly's theorem. Since Helly's theorem implies 
Caratheo\-dory's representation
theorem (Eggleston \cite{eggleston}), Proposition 2 in \S 
5 extends
also Caratheodory's theorem to the same wide range of 
families. Corollary 4
is the Knaster-Kuratowski-Marzukiewicz theorem (Berge 
\cite[p.\ 173]{berge},
which follows immediately
as a very special case of our Theorem 7 in \S 5. KKM's 
theorem
is restricted to families of sets in $R^n$ which cover an 
$n$-simplex,
while our Theorem 7 applies to families in a general 
topological space of
any cardinality, which cover any simplicial complex. An 
additional extension
of the KKM is Corollary 5, which applies to simple 
families. Corollary 6 is
the Brouwer fixed point theorem (Hirsch \cite{hirsch}). 
These results exhibit
a common topological root for these classical and somewhat 
disparate results.

\begin{corollary}[Helly's theorem]
Let $\{U_\alpha \}_{\alpha \epsilon S}$ be a family of
convex sets in $R^n$ with at least $n+1$ elements. Then if 
every subfamily
with $n+1$ sets has a nonempty intersection, the whole 
family has a
nonempty intersection.
\end{corollary}

\begin{pf} This follows directly from Proposition 2 in
\S 5, which is valid in much more generality 
for any number of sets in a general topological
space, because convex sets define an excisive family. 
%%$\diamondsuit $
\end{pf}

The following corollary requires no convexity:

\begin{corollary}[KKM Theorem]
Let $\{C_\alpha \}_{\alpha \in S}$ be a regular cover of a 
$k$-%
simplex $X$ as defined in \S \RM2. Then $\bigcap_{\alpha 
\in S}C_\alpha
\neq \emptyset .$
\end{corollary}

\begin{pf} This follows directly from Theorem 8, which is
valid more generally for any simplicial complex. Since 
nerve $\{C_\alpha
\}_{\alpha \in S}$ and $X$ have the same combinatorial 
structure, it
follows, in particular, that $\bigcap_{\alpha \in 
S}C_\alpha \neq \emptyset $.
%.\diamondsuit $
\end{pf}

In addition, the following result extends the KKM theorem 
to a different
class of covers, {\em simple covers}, as defined in \S 2, 
which need not
satisfy any of the conditions of KKM theorem:

\begin{corollary}[Extension of KKM to simple families]
 Let $\{U_\alpha \}_{\alpha \in S}$ be
a simple cover of a $k$-dimensional simplex $X$. Then 
$\bigcap_{\alpha \in
S}U_\alpha $ is not empty.
\end{corollary}

\begin{pf} This follows directly from Theorem 7, which is
also valid for covers of any complex $X$. %$\diamondsuit $
\end{pf}

Since the KKM theorem follows directly from Theorem 8 as 
shown in
Corollary 4, by presenting for completeness a well-known 
argument, we show
that Brouwer's fixed point theorem also follows as an 
immediate corollary of
our Theorem 8.

\begin{corollary}[Brouwer's fixed point theorem]
 Let $X$ be a $k$-simplex, and $%
f:X\rightarrow X$ a continuous function. Then $\exists 
x\in X:f(x)=x.$
\end{corollary}

\begin{pf} The proof follows by contradiction. If 
$f:X\rightarrow X$ has no
fixed point, then it defines a retraction $r:X\rightarrow 
\partial X$. Let $%
\partial X=\bigcup_iX_i,$ where $X_i$ is the $i$\<th face 
of $X,$ a $k-1$
simplex. Now define the closed sets $C_i=\{r^{-1}(X_i)\},$ 
$i=1,\ldots, k+1.$ Then 
$\{C_i\}_{i=1,\ldots, k+1}$ is a closed cover of $X$ 
satisfying the conditions of
Corollary 4, so $\bigcap_iC_i\neq \emptyset .$ But if $p\in
\bigcap_iC_i $, then $r(p)\in \bigcap X_i=\emptyset ,$ a 
contradiction. 
%$ \diamondsuit $
\end{pf}

\section{Market equilibrium and social choice}

Our final task is to establish the equivalence of the two 
economic problems,
namely, 
the existence of a competitive equilibrium and the 
existence of a social
choice function. A good way to start is to provide 
examples of spaces of
preferences in order to illustrate the topological problem 
involved in
social choice. By Theorem 3 in \S 4, this problem can be 
solved only 
for acyclic spaces of preferences.

A{\em \ preference} $\rho $ is an ordering of the choice 
space $R^n$ which
is induced by a utility function $u:R^n\rightarrow R,$ 
where we indicate $%
x\succeq _\rho y \Leftrightarrow u(x)\geq u(y).$ A {\em 
smooth preference}
on $R^n$ is defined by a smooth $(C^2)$ unit vector field 
$\rho
:R^n\rightarrow S^{n-1},$ with the property that $\exists 
$ a function $%
u:R^n\rightarrow R$ such that $\forall x\in R^n,$ $\exists 
\lambda (x)>0$
such that $\rho (x)=\lambda (x)Du(x)$; i.e., there exists 
a function $u$ such
that $\forall x,$ $\rho (x)$ is collinear with the 
gradient of $u$ (see
Debreu \cite{debreu}).

One example of a space of preferences $P$
is the{\em \ space of all smooth preferences }on 
$R^n,$denoted $\Gamma
(R^n), $ endowed with the sup norm, $\Vert \rho -\kappa 
\Vert = 
\sup _{x\in
R^N}\Vert \rho (x)-\kappa (x)\Vert $. Another example of a 
space of
preferences is the space $P_L$ of {\em all} {\em linear 
preferences }on{\em \ }$%
R^n $, which are those preferences induced by linear 
utility functions on $%
R^n,$ $n>2$. The space $P_{L\text{ }}$ is the sphere 
$S^{n-1}.$ If the zero
preference is also included, we have the space $P_{LN}$ of 
all linear
preferences on $R^n$\<---this space is $S^{n-1}\bigcup 
\{0\}.$ Different
preference spaces arise in different applications (for 
examples, see, e.g.,
Heal \cite{heal}). Typically, preference spaces are not 
linear nor convex or
acyclic; for example, the space of smooth preferences 
$\Gamma (R^n)$ is not
acyclic \cite{chichilnisky}.

Our last task is to establish the connection between the 
existence of a
market equilibrium and the existence of a social choice 
function. Both
problems depend on the characteristics of the traders' 
preferences, but they
do so in two different ways. The market $E$ has a finite 
set of preferences,
one for each trader, $\{\rho _1,\ldots, \rho _H\}.$ The 
set of preferences in
the economy is therefore a discrete finite set of points 
in the space of
smooth preferences $\Gamma (R^n)$ defined above$.$ The 
social choice
function, by contrast, is generally defined on large 
spaces describing a
universe of all possible preferences, typically a 
connected subset $P$ of
the space of all smooth preferences in $\Gamma (R^n),$ 
which is not a finite
set.

In order to exhibit the connection between the two 
problems---the existence
of market equilibrium and that of a social choice 
function---we define a
space consisting of preferences which are naturally 
``close'' to those of
the preferences of the traders in the economy $E.$ The 
space of preferences $%
P_E$ consists of a large number of preferences, assumed to 
be
a connected subspace
of $\Gamma (R^n)$, all of which are, in a well-defined 
sense, similar to the
preferences in the market $E$. We therefore need to define 
what is meant by
a smooth preference which is similar to the preferences of 
the traders in
the market $E.$

A smooth preference $\rho \in P$ defined over allocations 
in $R^{nH}$ is
called {\em similar} to the preference of trader $i\in E$ 
in position $j$
when $\forall x\in R^{nH}$, the projection of $\rho (x)$ 
on the $j$\<th copy
of $R^n$ is in the market cone of trader $i$; i.e., 
$\forall x\in R^{nH},$ $%
\rho ^j(x)\in D_i.$ The interpretation of this condition 
is that the
preference $\rho $ increases in the direction of that of 
the trader $i$ in
position $j $ for large utility values. The space 
$P(E_\theta )$ of
preferences similar to those of a subset $\theta \subset 
\{1,\ldots, H\}$ of
traders in $E$ was already defined in \S 4; it consists of 
all those
smooth preferences $\rho \in \Gamma (R^{nH})$ such that 
$\forall x\in
R^{nH}, $ $\rho ^j(x)\in \bigcup_{i\in \theta }D_i.$ If we 
consider the
problem of finding a social choice function for the space 
of preferences $%
P(E_\theta )$ which are similar to those of some subset 
$\theta $ of traders
in the economy $E,$ $\theta \subset \{1,\ldots, H\},$ then 
by Theorem 4 in \S
4 the necessary and sufficient condition is the acyclicity 
of $%
\bigcup_{i\in \theta }D_i.$ The existence of a social 
choice function for
every such space of preferences $P(E_\theta ),$ $\forall 
\theta \subset
\{1,\ldots, H\}$ therefore requires 
$$
\forall \theta \subset \{1,\ldots, H\}, \quad \theta \neq 
\emptyset \Rightarrow
\bigcup_{i\in \theta }D_i\text{ is acyclic.} 
$$

Note that in order to solve the social choice problem we 
must go back to the
properties of the family of market cones $\{D_i\}_{i\in 
\{1,\ldots, H\}}$ of the
economy $E$ defined in (\ref{t})---the same family of 
cones which define the
condition of limited arbitrage (\ref{la}).

Theorem 11 exhibits the identity between the problems of 
existence of a
competitive equilibrium for a market $E$ and the existence 
of a social
choice function. Let $E$ be a market as defined in \S 3. A 
{\em %
subeconomy} $E_\theta $ of $E$ is the market consisting of 
the those traders
in $E$ who belong to the set $\theta \subset \{1,\ldots, 
H\}$, i.e.,%
$$
E_\theta =\{R^n,\Omega _i,u_i,i\in \theta \}. 
$$

\begin{theorem}
The following properties of the economy $E=\{R^n,\Omega 
_i,\rho _i,i=1,\ldots, H\}$
are equivalent\/\RM: 

\RM{(a)} $E$ has a competitive equilibrium.

\RM{(b)} Every subeconomy $E_\theta$ of $E$ has a 
competitive equilibrium.

\RM{(c)} Every subeconomy $E_\theta$ of $E$ with at
most $n+1$ traders has a competitive equilibrium.

\RM{ (d)} There exists a social
choice function $\Phi :P(E_\theta )^k\rightarrow 
P(E_\theta )$ satisfying
conditions \RM{4.2.1} and \RM{4.2.2}, for every space 
$P(E_\theta )$ of
preferences similar to those of the traders in a nonempty 
set $\theta ,$ $%
\forall \theta $$\subset \{1,\ldots, H\},$ and $\forall 
k\geq 1.$
\end{theorem}

\begin{pf} The equivalence between (a) and (b) follows
immediately from Theorem 1 in \S 3 and from the definition 
of limited
arbitrage $(\LA)$ in (\ref{la}). We establish next the 
equivalence of the
statements (a) and (c). By Theorem 1, $E$ has a 
competitive equilibrium if
and only if $E$ has limited arbitrage $(\LA)$ as defined 
in (\ref{la}), i.e., 
if and only if the family of dual cones 
$\{D_i\}_{i=1,\ldots, H}$ has a nonempty
intersection. Since $\{D_i\}_{i=1,\ldots, H}$ is an 
acyclic excisive 
family in $%
R^n $, by Corollary 1, (\ref{la}) is true if and only if 
every subfamily of $%
\{D_i\}_{i=1,\ldots, H}$ with indices in a set $\theta 
\subset \{1,\ldots, H\}$ of at
most $n+1$ elements has nonempty intersection, i.e., if 
and only if the
corresponding subeconomy $E_\theta $ satisfies limited 
arbitrage (\ref{la}),
and therefore by Theorem 1 
if and only if $E_\theta $ has a competitive equilibrium.

The equivalence between statements (a) and (d) follows 
from Theorem 4 in
\S 4 and from Theorem 6 and Corollary 2 in \S 5, because $%
\{D_i\}_{i=1,\ldots, H}$ is an acyclic excisive family, so 
$$
\bigcap_{i=1}^HD_i\neq \emptyset \Leftrightarrow \forall 
\text{ nonempty }%
\theta \subset \{1,\ldots, H\},\text{ }\bigcup_{i\in 
\theta }
D_i\text{ is acyclic.}%
%\diamondsuit 
\qed
$$
\renewcommand{\qed}{}\end{pf}
%\newpage\ 

\section*{Acknowledgment}
I thank G. Heal, M. Hirsch, I. James, G. Segal, M. Shub 
and E. Spanier for
valuable comments.

\end{document}